\documentclass[12pt]{amsart}
\usepackage{amsmath, amssymb}
\usepackage{amsthm}
\usepackage{latexsym}
\usepackage{enumerate}
\usepackage{pstricks,pst-node,pst-coil}
\usepackage{hyperref}

\def\wwwlink#1{{\htmladdnormallink{\tt #1}{#1}}}
\def\mnote#1{}
\makeatletter
\@addtoreset{equation}{section}
\renewcommand{\theequation}{\the\c@section.\the\c@equation}
\renewcommand{\appendix}{\par
  \setcounter{section}{0}%
  \setcounter{subsection}{0}%
  \gdef\thesection{\@Alph\c@section}
  \gdef\theequation{\@Alph\c@section.\the\c@equation}
}
\makeatother
\newtheorem{theorem}[equation]{\bf T{\footnotesize HEOREM}}
\newtheorem*{theoremnonum}{\bf T{\footnotesize HEOREM}}
\newtheorem*{theooneone}{\bf T{\footnotesize HEOREM~\ref{theo.oneone}}}
\newtheorem*{theoremtwo}{\bf T{\footnotesize HEOREM~\ref{theo.two}}}
\newtheorem*{theosorth}{\bf T{\footnotesize HEOREM~\ref{theo.S.orth}}}

\newtheorem{proposition}[equation]{\bf P{\footnotesize ROPOSITION}}
\newtheorem{lemma}[equation]{\bf L{\footnotesize EMMA}}
\newtheorem{corollary}[equation]{\bf C{\footnotesize OROLLARY}}

\theoremstyle{definition}
\newtheorem*{definition}{\bf Definition}
\newtheorem*{remark}{\bf Remark}
\newtheorem*{notation}{\bf Notation}
\def\proof#1{{\par\medbreak\noindent {\bf Proof\setbox0\hbox{#1}%
\ifdim\wd0=0pt .\else\ \ignorespaces #1.\fi}\enspace}}
\def\iop#1{{\par\medbreak\noindent {\bf Idea of proof\setbox0\hbox{#1}%
\ifdim\wd0=0pt .\else\ \ignorespaces #1.\fi}\enspace}}

\newtheorem*{example}{Example}
\newtheorem*{examples}{Examples}

\def\la{\lambda}
\let\de\delta
\let\ep\varepsilon
\let\epsilon\varepsilon
\let\phi\varphi
\let\om\omega
\let\si\sigma
\def\Si{\Sigma}

\def\qed{{\leavevmode\unskip\nobreak\hfil\penalty 50\hskip 1em%
  \hbox{}\nobreak\hfil\lower 1pt\hbox{$\Box$\kern-.5pt}\parfillskip 0pt
  \finalhyphendemerits 0\par\bigbreak}}
\def\res#1#2{{#1}\lower .11ex\hbox{$|$}\lower .644ex\hbox{$\scriptstyle
#2$}}
\DeclareMathAlphabet{\doba}{U}{msb}{m}{n}
\def\mN{\doba N}
\def\mR{\doba R}
\def\mC{\doba C}
\def\mZ{\doba Z}

\def\cS{{\mathcal S}}
\def\cE{\mathcal{E}}
\def\cF{\mathcal{F}}

\def\eref#1{(\ref{#1})}
\let\ti\tilde
\let\witi\widetilde

\def\Spin{\mathop{\rm Spin}}
\def\SO{\mathop{\rm SO}}
\def\SU{\mathop{\rm SU}}
\def\U{\mathop{\rm U}}
\def\GL{\mathop{\rm GL}}
\def\Pspin{P_{\rm Spin}}
\def\PSO{P_{\rm SO}}
\def\End{\mathop{\rm End}}
\def\diam{\mathop{\rm diam}}
\def\Stab{{\mathop{\rm Stab}}}

\newcommand {\Ric}{\operatorname{Ric}}
\newcommand {\Scal}{\operatorname{scal}}
\let\scal\Scal
\newcommand {\Sec}{\operatorname{sec}}
\let\sec\Sec

\newcommand {\Vol}{\operatorname{vol}}
\newcommand {\dVol}{\operatorname{dvol}}
\newcommand {\area}{\operatorname{area}}

\newcommand {\divv}{\operatorname{div}}
\newcommand {\Id}{\operatorname{Id}}
\newcommand {\id}{\operatorname{id}}
\newcommand {\rank}{\operatorname{rank}}
\newcommand\dirac{{\not \!\!D}}

\let\<\langle
\let\>\rangle
\let\na\nabla

\newcommand{\vb}{V}
\newcommand{\subspace}{W}

\begin{document}

\title{Manifolds with small Dirac eigenvalues are nilmanifolds}

\author{Bernd Ammann}
\address{{\hskip-\parindent
        Bernd Ammann\\
        MSRI Berkeley\\
        17, Gauss Way\\
        Berkeley, CA 94720-5070\\
        USA\\}\newline {\tt {www.berndammann.de/publications}}}

\author{Chad Sprouse}
\address{\hskip-\parindent
        Chad Sprouse\\
        Department of Mathematics\\
        CSU Northridge\\
        18111 Nordhoff Street\\
        Northridge, CA 91330-8313, USA}
\email{chad.sprouse@csun.edu}

\thanks{Research at MSRI is supported in part by NSF grant DMS-9810361.}

%\begin{center}
%  \datverp
%\end{center}

\begin{abstract}
Consider the class of $n$-dimensional Riemannian spin manifolds with bounded
sectional curvatures and diameter, and almost non-negative scalar curvature.
Let $r=1$ if $n=2,3$ and $r=2^{[n/2]-1}+1$ if $n\geq 4$.
We show that if the square of the Dirac operator on such a manifold
has $r$ small eigenvalues, then the manifold is
diffeomorphic to a nilmanifold and has trivial spin structure.
Equivalently,
if $M$ is not a nilmanifold
or if $M$ is a nilmanifold with a non-trivial spin structure,
then there exists a uniform lower bound on the $r$-th
eigenvalue of the square of the Dirac operator.
If a manifold with almost nonnegative scalar curvature has one small Dirac eigenvalue,
and if the volume is not too small, then we show that the
metric is close to a Ricci-flat metric on $M$ with a parallel spinor.
In dimension $4$ this implies that $M$ is either a torus or a
$K3$-surface.
\end{abstract}

\maketitle


\tableofcontents

%%%%%%%%%%%%%%%%%%%%%%%%%%%%%
\section{Introduction}
%%%%%%%%%%%%%%%%%%%%%%%%%%%%%

%v2.3 B modified
The theorem of Bochner implies that if a connected compact Riemannian manifold has
non-negative curvature operator,
then $b^p(M) \leq {n\choose p}$. Furthermore if $b^p(M)={n\choose p}$
for some $p$ between $0$ and $n$
then $M$ is isometric to a flat torus.
In \cite{petersen.sprouse:01} it was shown that Riemannian manifolds with
$almost$-nonnegative curvature operator, and ${n\choose p}$
eigen-$p$-forms with small eigenvalue
must be diffeomorphic to a nilmanifold. Given the Hodge-de Rham theorem, this could be
viewed as a quantitative generalization of Bochner's theorem.

Here we discuss a similar result for the Dirac operator on Riemannian spin manifolds.
Let $\lambda_i(\dirac^2)$ denote the $i$-th eigenvalue of the square of the
Dirac operator, and let
$\lambda_i(\nabla^*\nabla)$  denote the $i$-th eigenvalue of the connection
Laplacian on spinors.
Here and throughout the article we assume that all eigenvalues are counted with
multiplicity.
%v2.3
All manifolds are connected.
%v2.3
Let $r(n)=2^{[{n\over 2}]-1}+1$ for $n\geq 4$ and $r(n)=1$ for $n\leq 3$.
Our main result is:

\def\contentoneone{Let $(M^n,g, \chi)$ be a compact Riemannian spin manifold with
$|\Sec|<K$, $\diam < D$. Then there
is $\varepsilon=\varepsilon(n,K,D)>0$,
such that if  $\lambda_r(\nabla^*\nabla)<\varepsilon$,
then $M$ is diffeomorphic to a nilmanifold.
Furthermore, $\chi$ is the trivial spin structure on $M$.
}

\begin{theorem}\label{theo.oneone}
\contentoneone
\end{theorem}

Using the Schr\"odinger-Lichenerowicz formula $\dirac^2=\nabla^*\nabla+
\scal/4$ this implies:

\begin{corollary}\label{cor.dir}
Let $(M^n,g, \chi)$ be a compact Riemannian spin manifold with $|\Sec|<K$, $\diam < D$.
Then there
is $\varepsilon=\varepsilon(n,K,D)>0$,
such that if $\scal > -\varepsilon$ and $\lambda_r(\dirac^2)<\varepsilon$
then $M$ is diffeomorphic to a nilmanifold. Furthermore, $\chi$ is the trivial spin structure on $M$.
\end{corollary}

By reformulating we obtain.
\begin{corollary}\label{cor.low}
If $(M, \chi)$ is not spin-diffeomorphic to a nilmanifold
with a trivial spin structure,
then among all metrics with bounded diameter and curvature,
there is a uniform lower bound on the $r$-th $\dirac^2$-eigenvalue and $r$-th $\nabla^*\nabla$-eigenvalue
  $$\la_r(\nabla^*\nabla)\geq \ep=\ep(n,\max |\sec|,\diam)>0.$$
In particular, any metric $g$ on $M$ with $\scal > -4\ep(n,\max|\sec|,\diam)$
has a non-trivial uniform lower bound on the $r$-th Dirac eigenvalue
  $$\la_r(\dirac^2)\geq \ep(n,\max |\sec|,\diam)+ {\min \scal\over 4}.$$
\end{corollary}

Recall that the Atiyah-Singer index theorem implies
$\dim \ker \dirac\geq |\hat A(M)|$. We obtain.

\begin{corollary}
If $M$ is an $n$-dimensional compact spin manifold with $|\hat A (M)|\geq
r(n)$, then $M$ does not carry a metric with $\scal >
-4\ep(n,\max|\sec|,\diam)$ for the above $\ep>0$.
\end{corollary}

We will give some examples that explain the special role of nilmanifolds and why
we cannot replace $r$ by a smaller number.
\begin{examples}\ \\[-5mm]
\begin{enumerate}[{\rm (1)}]
%\refstepcounter{enumi}
%\setcounter{enumi}{4}
\item \label{exam.nil}
%v2.3
Any nilmanifold $M^n$
carries a sequence of ``left-invariant'' metrics
$g_i$ with $\max |\sec_i| \to 0$, $\diam_i\to 0$ and $\Vol_i\to 0$.
If the spin structure on $M$
is trivial, then for this sequence of metrics $\la_{s}(\dirac_i^2)\to 0$ where
$s = \rank (\Si_i M)= 2^{[n/2]}>r$.
\item \label{kthreeex}
Let $N$ be a K3-surface. For an integer $n\geq 4$, let
$M$ be the Riemannian product $N\times T^{n-4}$, where the $n-4$ dimensional
torus carries an arbitrary flat metric. We equip $M$ with the product spin structure
of the unique spin structure on $N$ and the trivial spin structure
on $T^{n-4}$. Then the spinor bundle $\Si M$  on $M$ is isomorphic
(as a metric bundle with connection) to
$\pi_1^*(\Si N) \otimes_\mC \mC^{\tilde r}$ where
$\tilde r= 2^{[(n-4)/2]}=2^{[n/2]-2}$ and $\pi_1:N\times  T^{n-4}\to N$ is the
projection to the first component. If $\psi\in\Gamma(\Si N)$
is a parallel spinor on $N$ and $v$ is a constant section of $\mC^{\tilde r}$,
then $\pi_1^*(\psi)\otimes v$ is a parallel spinor on $M$.
As $N$ carries a 2-dimensional space of parallel spinors, the dimension of
the space of parallel spinors on $M$ is at least
$2\tilde r= r-1$. And hence
  $$\lambda_1(\dirac^2)=\ldots = \lambda_{r-1}(\dirac^2)=0.$$
This example shows that we cannot replace $r$ by $r-1$ in the above theorem.
\end{enumerate}
\end{examples}
The next example will show that the sectional curvature bounds in
Theorem~\ref{theo.oneone} are necessary. We need a lemma.
\def\BR#1{{B^{{\mathbb R}^4}_{#1}(0)}}
\begin{lemma}\label{lem.cyl}
Let $h$ be the standard metric on $S^3$.
The ball of radius $R$ around~$0$ in Euclidean space ${\mathbb R}^4$
will be
denoted as $\BR{R}$.
For any $\epsilon,R,\rho>0$ there is $a>0$ and a metric
$g=dt^2 + \phi^2(t) h$ on
$S^3\times (-(a+\rho),(a+\rho))$ such that
\begin{enumerate}[{\rm (a)}]
\item $g|_{(-(a+\rho),-a]\times S^3}$ and $g|_{[a,(a+\rho))\times S^3}$
are isometric to $\BR{R+\rho}\setminus \BR{R}$,
\item ${\rm scal}_g\geq -\epsilon$, and
\item ${\rm diam }(M,g)\leq  6(R+\rho)$
\end{enumerate}
\end{lemma}

For the proof one translates the desired properties into an ordinary
differential inequality for $\phi$. % which is easy to solve.
Details are available in \cite{ammann.sprouse:detailssmallev1:1}.

Using the lemma we can construct an example showing that the curvature bound is necessary.

\pagebreak[3]
\begin{example}\ \\[-5mm]
\begin{enumerate}[{\rm (1)}]
%\refstepcounter{enumi}
\setcounter{enumi}{2}
\item
Consider the flat torus $T^4$. Let  $Z:=Z(R_j,\ep_j):=(-c,c)\times S^{n-1}$
carry a metric as in the above lemma
%in Lemma~\ref{lem.cyl}
%v2.3
for $\ep_j$ and $R_j$ sufficiently small, that we will choose later, and for $\rho_j=R_j$.
Let $(M_j,g_j)$ be given by removing $2j$ small disks from $T^4$
and attaching $j$ handles isometric to $Z$.
% ($\ep$ and $R$ depend on $j$).
The trivial spin structure on $T^4$
can be extended to a spin structure on $M_j$.\footnote{This extension is not
unique, but our construction works for any choice of spin structure.}
For a suitable choice of $R_j$ and $\ep_j$
we obtain a family of Riemannian manifolds $(M_j,g_j)$ with uniformly bounded
diameter and $\liminf_{j\to \infty} \min \scal_j=0$.
They are pairwise non-diffeomorphic, and the sectional
curvature is not uniformly bounded.
Following the lines of \cite{baer.dahl:02} we use a cut-off function
vanishing in the handles to construct a $(\rank \Si M_j=4)$-dimensional space
of test spinors.
From this we see that $D^2$ has at least~$4$
eigenvalues arbitrarily close to~$0$.
This example shows that the sectional bound in Theorem~\ref{theo.oneone}
is necessary if the dimension of $M$ is $4$. We obtain similar examples
for higher dimensions $n$ by taking the product with $(S^1)^{n-4}$.
See
%Appendix~\ref{appendix.cyl} and
\cite{ammann.sprouse:detailssmallev1:1}
for details.
\end{enumerate}
\end{example}
Finally, we will give several examples in order to show that this bound generalizes
previously known bounds in dimension $2$ and $3$.

\begin{examples}\ \\[-5mm]
\begin{enumerate}[{\rm (1)}]
%\refstepcounter{enumi}
\setcounter{enumi}{3}
\item If $M$ is diffeomorphic to the 2-dimensional sphere,
then such a lower bound is already known.
It is a result of B\"ar~\cite{baer:92} that
  $$\la_1(\dirac^2)\geq {4\pi\over \area(M)}.$$
If $K\geq -\de^2$, $\de> 0$ then
$\area(M)\leq {2\pi\over \de^2}\,\left[\cosh\left({\de \diam M}\right)-1\right]$.
%$\area(M)\leq {2\pi\over \de^2}\,\left[\cosh\left({\de \diam M\over 2}\right)
%-1\right]$.
%%$\pi \diam (M)^2/4+ \tau(\de| D)$ where $\tau$ is a
%%continuous function with
%%$\tau(\de| D)\to 0$ for $\de\to 0$.
Hence,
  $$\la_1(\dirac^2)\geq {2\de^2\over \cosh\left({\de \diam M}\right)-1}
    ={4 \over  \diam M^2}- O(\de^2\diam M^2).$$
\item If $(M,g)$ is diffeomorphic to the 2-dimensional torus $T^2$ equipped with a non-trivial spin structure,
then it is not difficult to derive an explicit lower bound on $\la_1(\dirac^2)$
from previously known estimates. To derive this,
we use the uniformization
theorem to find $u\in C^\infty(T^2)$ and a flat metric
$g_0$ with $g=e^{2u}g_0$. Using the estimates in \cite{ammann:02} together with
some elementary calculations\footnote{Details of this calculation available on
\wwwlink{http://www.berndammann.de/publications/smallev1}.} one obtains
  $$\max u -\min u \leq \cS(K,D),$$
where $\cS$ is an explicitely known, but long expression with $\cS(0,D)=0$. Then one
easily derives from \cite[Corollary~2.3]{ammann:02} that
  $$\la_1(\dirac^2)\geq {\pi^2\over 4D^2}\,e^{-4\cS(K,D)}.$$
\item Surfaces of genus greater than 1 cannot have almost non-negative curvature
in the above sense. Hence, (1) and (2)
yield an explicit, but long formula for $\ep$ in dimension $n=2$.
However, in higher dimension one expects $\ep$ to be an even more complicated
expression. Thus, we want to restrict our attention to existence results.
\item Let $(M,g,\chi)$ be a compact spin $3$-manifold with $\scal\geq 0$
and $\ker \dirac\neq \{0\}$. Because of
the Schr\"odinger-Lichnerowicz formula, any $\phi\in \ker D\setminus
\{0\}$ is a nontrivial
parallel spinor, which implies that $(M,g)$ is Ricci-flat, and hence flat.
However, any flat compact $3$-manifold admitting a nontrivial parallel spinor is diffeomorphic to a torus
(see \cite[Theorem~5.1]{pfaeffle:00}) and the spin structure is the trivial one.
\end{enumerate}
\end{examples}

We compare Corollary 1.2 which gives a uniform lower $r$-th eigenvalue
bound to a theorem of J.~Lott which provides
a uniform upper bound on all eigenvalues.
\begin{theoremnonum}[{Lott \cite[Theorem 4]{lott:02}}]
Let $k\in \mZ^+$. Then there is an $E_k=E(n,K,D,k)$ such that any compact
Riemannian spin manifold $(M,g,\chi)$ with $|\sec|<K$, $\diam \leq D$ satisfies either
\begin{enumerate}[{\rm (a)}]
\item $\la_k(\dirac^2)\leq E_k$
\item $M$ is the total space of an affine
fiber bundle $M\to B$ with possible singularities, whose generic fiber
is an infranilmanifold, and the spin structure along the generic fibers is \textbf{not trivial}.
\end{enumerate}
\end{theoremnonum}

Another result which will be proven in section~\ref{sec.einstein}
gives a different conclusion for manifolds with only one small Dirac eigenvalue, and
additionally, a lower volume bound.

\begin{theorem}\label{theo.two}
%v2.3
Let $(M,g)$ have $|\Sec| < K$, $\diam < D$, $\Vol > v$.
Let $\lambda_1(\nabla^*\nabla)$ denote the first eigenvalue
of the connection Laplacian $\nabla^*\nabla$ on the spinor bundle with respect to a spin
structure $\chi$.
Then for all $\delta>0$, there is an $\varepsilon=\varepsilon(n,v,K,D,\delta)>0$ such that
if $\lambda_1(\nabla^*\nabla) < \varepsilon$, then $(M,g,\chi)$
has $C^{1,\alpha}$-distance $\leq \delta$ to a Ricci-flat Einstein metric with a
nontrivial parallel spinor.
\end{theorem}

\begin{corollary}\label{cor.two}
%v2.3
For $\delta >0$, there is an $\varepsilon=\varepsilon(n,v,K,D,\delta)>0$
such that the following holds:
Let $(M,g)$ have $|\Sec| < K$, $\diam < D$, $\Vol > v$ and
$\Scal > -\epsilon$.
Let $\lambda_1(\dirac^2)$ denote the first eigenvalue of $\dirac^2$ with
respect to a spin structure $\chi$.
If $\lambda_1(\dirac^2) < \varepsilon$, then $(M,g,\chi)$ has $C^{1,\alpha}$-distance
$\leq \delta$ to a Ricci-flat Einstein metric with a nontrivial parallel spinor.
\end{corollary}

A compact $4$-dimensional manifolds $M$ carrying a parallel spinor is either
a flat torus or a $K3$-surface. Hence, any $4$-manifold with one
small Dirac eigenvalue is either diffeomorphic to a torus or a
K3-surface, or is collapsed.

Example~\eref{exam.nil} shows that the volume bound in the above theorem and corollary
is necessary.

The structure of the article is as follows.
In Section~\ref{sec.gen} we will reformulate some previously
known estimates on vector bundles. In Section~\ref{sec.einstein} we will
apply these estimates to prove Theorem~\ref{theo.two}.
In the following sections Theorem~\ref{theo.oneone} is proved.
We will develop most of the tools in such a generality
that we can easily replace the Dirac operator (acting on spinors)
by other elliptic operators acting on sections of bundles with special
holonomy. We begin this in Section~\ref{sec.fix.faith} by defining the
fixing dimension $r$ of a faithful representation, which immediately gives
the fixing dimension $r$ of a vector bundle with special holonomy.
In Section~\ref{sec.vb.hol} we show that if there are $r$ almost
parallel sections on such a bundle, then the bundle is trivialized by almost
parallel sections. Section~\ref{sec.spin.fix} determines the fixing number
of the spinor bundle, and we finally prove Theorem~\ref{theo.oneone} in
the last section.

%Note that with similar techniques one can also classify manifolds
% with positive
%scalar curvature and sufficiently many eigenvalues of the Dirac operator
%that are sufficiently close to T.~Friedrich's Dirac eigenvalue bound.
If one applies the techniques that we will present in this paper to the
Friedrich connection on the spinor bundle instead of the standard connection,
one obtains analogs of
Theorem~\ref{theo.oneone}, Corollary~\ref{cor.dir}, Theorem~\ref{theo.two}
and Corollary~\ref{cor.two}. In particular, we obtain the following theorem.
Here, once again, we define $r=1$ if
$n=2,3$ and $r=2^{[n/2]-1}+1$ if $n\geq 4$.

\begin{theorem}
Let $(M^n,g,\chi)$ be a compact Riemannian manifold
with $\diam < D$, $| \Sec | < K$, and $\scal \geq n(n-1)\rho^2$
with a constant $\rho > 0 $. Let $\dirac$ be the Dirac operator on $M$.
Then for any $\delta>0$, there is $\varepsilon = \varepsilon(n,K,D,\rho,\delta)$
such that if $\dirac$ has $r=r(n)$
eigenvalues $\lambda_i\in [0, \frac{n\rho}{2} + \varepsilon)$,
then $M$ has $C^{1,\alpha}$-distance $\leq \delta$ to a manifold of
constant curvature $\Sec \equiv \rho^2$.
\end{theorem}

However, motivated by B\"ar's classification of manifolds with real
Killing spinors \cite{baer:93}, we conjecture that the theorem still holds
for a smaller number $r$. This is ongoing research.

Again, one sees that the bound on the curvature is necessary.
\begin{example}\ \\[-5mm]
\begin{enumerate}[{\rm (1)}]
%\refstepcounter{enumi}
\setcounter{enumi}{7}
\item On any compact manifold $M$ that admits a metric of positive
scalar curvature, and an arbitrary spin structure on $M$
\mnote{B2C: I have to verify the details},
C.~B\"ar and M.~Dahl \cite{baer.dahl:p03} have constructed
a sequence of metrics $g_i$ on $M$ with scalar curvature $\geq n(n-1)$,
but with $\lambda_{2^{[n/2]}}(\dirac_{g_i}^2)\to n^2/4$.
\end{enumerate}
\end{example}

\emph{Throughout the paper we adopt the convention that
$\tau(x_1|x_2,\ldots,x_m)$ represents a continuous function
in $x_1,\ldots,x_m$ such that $\tau \to 0$ as $x_1 \to 0$ with $x_2,\ldots, x_m$
fixed.}

{\bf Acknowledgement.} We want to thank C.~B\"ar, A.~Degeratu,
M.~Dahl, P.~Petersen,
and W.~Tuschmann for several helpful discussions. The article was
completed while the first named author enjoyed the hospitality of
the MSRI Berkeley, CA, USA. Research at MSRI is supported in part by 
NSF grant DMS-9810361.

\section{General estimates on vector bundles}\label{sec.gen}
Let $\vb$ be a complex vector bundle of rank $k$
over $M$ equipped with a connection
$\nabla$ and a metric $\<\cdot,\cdot\>$.
Recall that the second covariant derivative on sections of $\vb$ is given by
$\nabla^2_{X,Y} S = \nabla_X \nabla_Y S - \nabla_{\nabla_X Y}S$, and the
curvature tensor on sections of $\vb$ is given by
$R^\vb(X,Y)S = \nabla^2_{X,Y}S
- \nabla^2_{Y,X}S$.

Furthermore, we consider the connection Laplacian on $\vb$, which
is given by
$$ \nabla^*\nabla S = -\sum_{i=1}^n \nabla^2_{e_i,e_i}S, $$
where $\{e_i\}$ is an orthonormal set of vectors at any point $p \in M$.
We say that $S$ is an eigensection of $\vb$ with eigenvalue $\lambda$
if $\nabla^*\nabla S = \lambda S$.
In this general situation, we recall the
eigenvalue pinching theorems from (\cite{petersen.sprouse:99},
\cite{petersen.sprouse:01}),
which characterize eigensections with $small$ eigenvalues.

\begin{notation} We use the volume-normalized $L^p$-norm given by
  $$\|u\|_p = \left( \frac{1}{\Vol M} \int_M |u|^p \dVol\right)^{1/p}$$
and the volume-normalized $L^2$-scalar product
  $$(u,v)=\frac{1}{\Vol M} \int \bar u v\,\dVol.$$
\end{notation}

\begin{theorem}\label{theo.S.est}
Suppose that $\vb$ satisfies $|R^\vb|, |\nabla R^\vb| < K$, and $M$ satisfies $|\sec |< K, \diam < D$.
Suppose $S$ is an eigensection of $\vb$ with eigenvalue $\lambda$, normalized so that $\|S\|_2=1$. Then,
$$ \|S\|_\infty \leq 1 + \tau(\lambda|n,K,D)$$
$$\| \nabla S\|_\infty \leq \tau(\lambda|n,K,D)$$
$$ \| \nabla \nabla S\|_2 \leq \tau(\lambda|n,K,D)$$
\end{theorem}
And furthermore,
\begin{theorem}\label{theo.S.orth}
%\contenttheosorth
Suppose that $S_1,\dots,S_m$ are $L^2$-orthonormal eigensections of $\vb$, with
eigenvalues $\lambda_1 \leq \dots\leq \lambda_m$.
Then with $K,D$ as above,
$$ \| \langle S_i, S_j \rangle - \delta_{ij}\|_\infty \leq \tau(\lambda_m|n,K,D).$$
\end{theorem}

We outline the proofs of these facts in Appendix A.

%%%%%%%%%%%%%%%%%%%%%%%%%%%%%
\section{The first Dirac eigenvalue}\label{sec.einstein}
%%%%%%%%%%%%%%%%%%%%%%%%%%%%%

In this section we characterize manifolds with $|\Sec|<K$, $\diam < D$, almost
non-negative scalar curvature, and a single small Dirac eigenvalue. Our result is that such manifolds are either collapsed in
the sense of Cheeger-Fukaya-Gromov or $C^{1,\alpha}$-close to an Einstein manifold
with a parallel spinor.
Note that in the case of the first eigenvalue on differential $p$-forms $\lambda_1^+(\Delta_p)$, such a manifold would always be collapsed.
That is, there is a lower bound on $\lambda_1^+(\Delta_p)$
given a lower volume bound and the above curvature and diameter bounds (\cite{CC:90}). This is proved as follows.
Suppose $(M_i,g_i)$ is a sequence of manifolds as above with $\lambda^p_1 \to 0$.
Then the above conditions imply that there is a subsequence of $M_i$ that converges
to a limit manifold $\overline M$
in the $C^{1,\alpha}$ topology. But this is not possible since $\overline M$ would have a
higher $p$-th Betti number than the limiting manifolds $M_i$.

On the other hand, since the number
of harmonic spinors is not topologically invariant,
this argument will not work in the
spinor case. The reader interested in harmonic spinors may consult
the classical reference \cite{hitchin:74} or several articles
containing recent results \cite{baer:96b},\cite{baer.dahl:02}, \cite{maier:97}
about the dependence of $\dim \ker \dirac$ on the metric.

\begin{theoremtwo}\mnote{B2Chad: I had to correct the role of the scalar
curvature here.}
Let $(M,g)$ have $|\Sec| < K$, $\diam < D$, $\Vol > v>0$.
Let $\lambda_1(\nabla^*\nabla)$ denote the first eigenvalue
of the connection Laplacian $\nabla^*\nabla$ with respect to any spin
structure.
Then for any $\delta>0$, there is $\varepsilon=\varepsilon(n,v,K,D,\delta)$ such that
if $\lambda_1(\nabla^*\nabla) < \varepsilon$, then $(M,g)$
has $C^{1,\alpha}$-distance $\leq \delta$ to a Ricci-flat Einstein metric with a
nontrivial parallel spinor.
\end{theoremtwo}

\proof{}
Using Theorem~\ref{theo.smoothing}
and Proposition~\ref{prop.conn.cont}
we can assume that $|\nabla R|<K_1(K,n,D)$.
Let $\sigma$ denote an eigenspinor to the eigenvalue
$\lambda_1$ with $\| \sigma \|_2=1$. Then
$\||\sigma|-1\|_\infty < \tau(\lambda_1|n,K,D)$ and
$|| \nabla \nabla \sigma ||_2 < \tau(\lambda_1|n,K,D)$.
If $e_1,\ldots,e_n$ denotes a local orthonormal frame, then from the curvature formula
for spinors
$R(X,e_i)\si= {1\over 4}\, \sum_{j,k}
\<R(X,e_i)e_j,e_k\>\,e_i\cdot e_j \cdot \sigma$ one deduces (see e.g.\ \cite{friedrich:80})
$$\Ric(X)\cdot\sigma = -2 \sum_{i=1}^n e^i \cdot R(X, e_i)\sigma,$$
and the fact that $R(\cdot, \cdot)\sigma$
is clearly bounded by $\nabla\nabla\sigma$,
this implies that $\| \Ric \|_2 < \tau(\lambda_1|n,K,D)$.
The lower bound on the volume together with the upper bounds
on $|\Sec|$ and $\diam$
imply a lower bound on the injectivity radius of $(M,g)$.
Then from \cite[Theorem 6.1]{petersen:97} (see also \cite{anderson:90},
\cite{gao:90}),
we have that for $\lambda_1 < \varepsilon(n,v,K,D)$, $M$ is
$C^{1,\alpha}$-close
to a $C^\infty$ Einstein manifold $\overline M$ with $\Ric \equiv 0$.
%v2.3
As $\overline M$ is diffeomorphic to $M$ if $M$ and $\overline M$
are $C^{1,\alpha}$-close, we may assume in the following that $M$ and
$\overline M$ are equipped with the same spin structure.

It remains to show, that if $\ep$ has been chosen small enough, then
$\overline M$ must carry a parallel spinor.
Assume the opposite,
then we have a sequence of manifolds $M_i$ converging to $\overline M$ in
the $C^{1,\alpha}$ topology, with $\lambda_1(\nabla^*\nabla,M_i) \to 0$.
%v2.3
Proposition~\ref{prop.conn.cont} implies
$\lambda_1(\nabla^*\nabla, \overline{M})=0$, or in other words
%This means that $\overline M$ has
%$\lambda_1(\nabla^*\nabla) < \tau(i^{-1})$ for all $i$ which contradicts
%the fact that $\nabla^*\nabla$ has discrete spectrum on $\overline M$,
%unless
$\overline M$ admits a parallel spinor.
\qed

\begin{remark}
The above proof can be slightly simplified by using spinors on
manifolds with a $C^{1,\alpha}$-metric. We avoided this for technical reasons.
\end{remark}

%%%%%%%%%%%%%%%%%%%%%%%%%%%%%
\section{The fixing dimension of a faithful representation}\label{sec.fix.faith}
%%%%%%%%%%%%%%%%%%%%%%%%%%%%%

Here we discuss the fixing dimension for representations which we will need in Section 5.
Let $G$ be a Lie group, and let $\rho:G\to \End(V)$ be a faithful (i.e. injective)
complex representation.
%Then for $r=1,\dots, \dim V$, let $\rho^k$ be the induced action on the
%\emph{$r$-frame bundle}
%  $$F_r:=\{(v_1,\dots,v_r)\,|\, v_1,\dots, v_r\mbox{ linearly independent}\}.$$
For any subspace $\subspace$ of $\vb$ let
  $$\Stab^G(\subspace):=\{g\in G\,|\, \rho(g) w=w\quad \forall w\in \subspace\}$$
be the stabilizer. Note that faithfulness of $\rho$ means that $\Stab^G(V)=\{1\}$.

\begin{definition}
The \emph{fixing dimension} $\cF(\rho)$ of $\rho$ is defined to
be the smallest number $r\in \{0,\dots,\dim V\}$
with the property that any $r$-dimensional subspace $\subspace\subset V$ has a finite stabilizer $\Stab^G(\subspace)$.
\end{definition}

For the standard representation of $U(n)$ on $\mC^n$, we have
  $$\cF\Bigl(U(n)\hookrightarrow \GL(\mC^n)\Bigr)=n,$$
whereas
 $$\cF\Bigl(SU(n)\hookrightarrow \GL(\mC^n)\Bigr)=n-1.$$
\begin{proposition}
Let $\rho$ be a unitary representation with fixing dimension $\cF(\rho)$.
Then there is $N\in \mN$ such that for all $\cF(\rho)$-dimensional $W\subset V$,
we have $\# \Stab^G(W)\leq N$. We denote by $N(\rho)$ the smallest such
$N$.
\end{proposition}

The Proposition follows from the following lemma
by letting $K$ be  the set of all orthonormal $k$-frames $(v_1,\dots,v_k)\in V^k$,
where $k=\cF(\rho)$.

\begin{lemma}
If a compact Lie group $G$ acts continuously on a compact manifold $K$,
such that for any $p\in K$ the stabilizer $\Stab^G(p)$ is finite,
then the $\# \Stab^G(p)$ is uniformly bounded on $K$, i.e.\  there is
$N=N(G,K)$ such that $\# \Stab^G(p)\leq N$ for all $p\in K$.
\end{lemma}

\proof{}
Suppose that there exists $p_i \in K$ such that $\# \Stab^G(p_i)\to \infty$. Then after choosing a subsequence we have
$p_i \to p$ for some $p \in K$. We will exhibit a $1$-parameter subgroup in $\Stab^G(p)$, which is hence an
infinite subgroup of $G$.
Since $\# \Stab^G(p_i)\to \infty$ we can choose
$g_i, \tilde g_i \in \Stab^G(p_i)$ such that $d(g_i, \tilde g_i) \to 0$ with respect to a left-invariant metric on $G$.
Then letting $h_i = g_i^{-1}\tilde g_i$ we also have that $h_i(p_i) = p_i$ and $d(e,h_i) \to 0$.
For each $i$, chose a unit-length $v_i \in \frak g$ and $t_i \in \mathbb R$
such that $\exp(t_iv_i)= h_i$. Hence $t_i\to 0$. Then again after choosing
a subsequence we can assume that $v_i \to v$ for some $v \in \frak g$. For any fixed $t \in \mathbb R$ choose a sequence of
integers $k_i$ with $k_i t_i \to t$. Then $\exp(k_i t_i v_i) \to \exp(tv)$. But $\exp(k_i t_i v_i)=h_i^{k_i}$, and
since $G$ acts continuously we have
that $(h_i^{k_i}, p_i) \to (\exp(tv), p)$ implies that $h_i^{k_i}p_i \to \exp(tv)p$. But $\lim_{i\to\infty} h_i^{k_i}p_i
= \lim_{i\to\infty}p_i = p$. Hence $\exp(tv)p=p$ for any $t\in \mathbb R$.
\qed

In section~\ref{sec.spin.fix} we will determine the fixing dimension of
the spinor representation.

%%%%%%%%%%%%%%%%%%%%%%%%%%%%%
\section{Eigenvalue pinching on vector bundles with special holonomy}\label{sec.vb.hol}
%%%%%%%%%%%%%%%%%%%%%%%%%%%%%

Let $\vb$ be a vector bundle of rank $k>0$
over $M$ equipped with a connection
$\nabla$ and a metric $\<\cdot,\cdot\>$.
We fix $p\in M$.
We assume that the holonomy group of the bundle is contained in a closed
Lie group $H\subset U(\vb_p)$. For any $q\in M$ let $H_q$ be
the parallel transport of $H$ to $q$.

\begin{proposition}\label{prop.hol.absch}
We assume that $M$ and $\vb$ satisfy the conditions of
Theorem~\ref{theo.S.est}.
Assume that the rank-$k$ bundle $\vb$ has holonomy contained in $H$.
Let $r=\cF(H\subset U(T_pM))$ be the fixing dimension of the holonomy.
Let $S_1,\dots,S_r$ be $L^2$-orthonormal sections of $\vb$ such that
  $$\na^*\na S_i = \lambda_i S_i,$$
$1\leq i \leq r$,
$0\leq \lambda_i\leq \epsilon$.
Then for small $\epsilon>0$, there is a  finite cover
$\pi:\widetilde {M}\to M$ and smooth sections
$e_1,\ldots,e_k$ of $\pi^*(V)$
with the following properties:
\begin{enumerate}[(1)]
\item $\cE:=(e_1,\ldots,e_k)$ is a frame, i.e.\ $\cE(q)$ is a basis of $\pi^*(\vb)_q$  for all $q\in \widetilde M$.
\item
  $$|S_i(\pi(q))-e_i(q)| \leq \tau(\epsilon|n,K,D)\qquad \forall q\in \widetilde M,
i=1,\dots,r .$$
\item
  $$(\na \cE)_q = (\na e_1,\dots,\na e_k)_q\in T^*_q\widetilde M\otimes \mathop{\rm Lie} (H_{\pi(q)})$$
\item
  $$|\na e_i(q)| \leq \tau(\epsilon|n,K,D)\qquad \forall i=1,\dots, k\quad\forall q\in \widetilde M .$$
\end{enumerate}

Furthermore, if
%$\Stab^{H}(\subspace)=\{1\}$ for all $r$-dimensional $\subspace$,
$N(H\subset U(T_pM))=1$, then we can choose $\widetilde M= M$.
\end{proposition}

\proof{}
According to Theorems~\ref{theo.S.est} and \ref{theo.S.orth} we have
the estimates
\begin{eqnarray}
   \|S_i\|_\infty & \leq & 1 + \tau(\ep|n,K,D)\label{S.est.rep}\\
   \| \nabla S_i\|_\infty & \leq & \tau(\ep|n,K,D)\label{naS.est.rep}
%\\
%   \| \nabla \nabla S_i\|_2 & \leq & \tau(\lambda|n,K,D)\|S_i\|_2\nonumber
\end{eqnarray}
and
\begin{equation}\label{S.orth.rep}
  \| \langle S_i, S_j \rangle - \delta_{ij} \|_\infty < \tau(\ep|n,K,D).
\end{equation}
We apply pointwise the Hilbert-Schmidt orthogonalization procedure
to $S_1,\dots,S_r$ and obtain new sections
$\tilde S_1,\ldots,\tilde S_r$.
All functions in this procedure and their first derivatives
are controlled in terms of $\ep$, $n$, $K$, and $D$.
As a consequence these new sections also satisfy
\eref{S.est.rep} and \eref{naS.est.rep} and are pointwise orthonormal.
We fix a point $p\in M$. The parallel transports of
$(\tilde S_1(p),\ldots,\tilde S_r(p))$
define a principal bundle over $M$ whose structure group
is the holonomy group.
By enlarging the structure group to $H$
we obtain an $H$-principal bundle which we will denote by $P_H(M)$.
%v2.3
The bundle $P_H(M)$ is a parallel subbundle of the frame bundle $P_{U(k)}(\vb)$.

We denote $G:=\Stab^{U(k)}(\mR^r)\cong U(k-r)$.
Let $P_G(M)$ be the bundle of orthonormal
bases of $\vb$ such that the first $r$ basis vectors coincide with
$\tilde S_1,\dots, \tilde S_r$ at each base point. Note that $P_G(M)$
is an $U(k-r)$ principal bundle. The bundles $P_G(M)$ and $P_H(M)$ have a
common point
$(\tilde S_1(p),\ldots,\tilde S_r(p))$.

Choose a bi-invariant metric on $U(k)$. This induces a
metric on each fiber of $P_{U(k)}(\vb)$. For $x\in M$ let $\delta(x)$
be the distance of the fiber of $P_G(M)$ over $x$ to the fiber of
$P_H(M)$ over $x$ with respect to this metric.
Then, $\delta:M\to [0,\infty)$ is a function with
$\delta(p)=0$. Using \eref{naS.est.rep} one sees that $\delta$ is a Lipschitz
function with Lipschitz constant of the form
$\tau(\ep| n,K,D)$.
%v2.3
We assume that $\ep$ is so small that $2\delta$ is
smaller than the injectivity radius of $U(k)$ and smaller than
$\inf \{d(A,e)\,|\,A\in G\cap H, A\neq e\}$.
Let $\widetilde M$ be the set of all elements of
$P_H(M)$ having minimal distance from $P_G(M)$. Because of symmetry
we have $\# G\cap H$ many points in $\widetilde M$
over each point in $M$.
%v2.3
As $\delta$ is chosen as above, $\widetilde M$
is a smooth manifold
and $\pi:\widetilde M\to M$ is a covering of $M$
with $\# G\cap H$ many leaves. Any $q\in\widetilde M\subset P_{U(k)}\vb$ can be
written as $q=(e_1(q),\dots,e_k(q))$ with $e_j\in V_{\pi(q)}$, and
$e_j$ are clearly smooth sections of $\pi^*(V)$ satisfying (1).

Because of our construction the distance between $e_i$ and $\ti S_i$
is bounded by $\de$, and hence we obtain (2).
As $\cE:=(e_1,\ldots,e_k)$ is a
section of $\pi^*(P_H(M))$, we see that (3) holds.

For (4) we have to prove that for $q$ in $M$
and an arclength-parametrized curve $c$ with $c(0)=q$,
\begin{eqnarray}\label{eq.four.to.prove}
  |\na_{\dot c(0)} e_i(q)|\leq \tau (\ep|n,K,D).
\end{eqnarray}
Let $\widehat \cE$ be the parallel transport of $\cE_q$ along $c$.
Obviously, $\widehat \cE\in P_H(M)$.
Suppose that $A\in U(n)$ is the unique matrix such that $\cE_q \cdot A$
is the closest point to $\cE_q$ in $(P_G(M))_q$. Then $\widehat \cE\cdot A$
is also a parallel frame of $\vb$ along $c$, and by the construction of
$P_G(M)$ and also \eref{naS.est.rep} one sees that the distance
from $\widehat \cE(t)\cdot A$ to $(P_G(M))_{c(t)}$ is bounded by
$|t|\tau (\ep|n,K,D)$.
By applying the implicit function theorem one concludes that the distance
between $\cE(c(t))$ and $\widehat \cE(t)$ is bounded by $|t|\tau (\ep|n,K,D)$
for small $t$. This implies \eref{eq.four.to.prove}.
\qed

\begin{figure}[htbp]
\begin{center}
\begin{pspicture}(0,0)(11,6)
%% Initializing
\newgray{grayone}{.75}
\newgray{graytwo}{.30}
\newgray{graythree}{.58}
%\psgrid(0,0)(11,6)

%%Curves
\pscurve(1,5)(3,4.6)(7,5.8)
\pscurve(1,0.5)(3.2,2.7)(7,3.2)
\psline(3.8,2.92)(3.6,4.7)
\psline(3.4,2.8)(3.0,4.4)
\psdot(3.0,4.4)

%%Right angles
\psarc(3.8,2.92){.3}{10}{100}
\psdot(3.9,3.05)
\psarc(3.6,4.7){.3}{280}{10}
\psdot(3.7,4.6)%%Punkte noch kleiner machen
\psarc(3.4,2.8){.3}{13}{103}
\psdot(3.5,2.94)

%%Notation
\rput(3.9,2.6){\psframebox*{$\cE$}}
\rput(3.7,5.1){\psframebox*{$\cE\cdot A$}}
\rput(3.5,2.5){${\widehat \cE}$}
\rput(2.4,4.3){\psframebox*{${\widehat \cE\cdot A}$}}
\rput(5.4,2.8){$(P_G(M))_{c(t)}$}
\rput(5.4,4.6){$(P_H(M))_{c(t)}$}

%%\rput(6.1,2.6){$T$}

\end{pspicture}
\end{center}

\caption{The shortest line between $(P_G(M))_{c(t)}$ and $(P_H(M))_{c(t)}$.}
\end{figure}
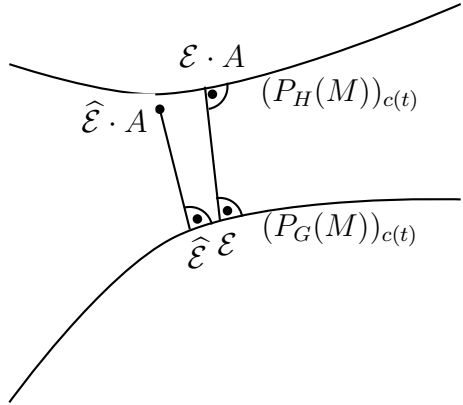

\section{The fixing dimension for the spinor representation}\label{sec.spin.fix}

Let $\vb$ be an $n$-dimensional real vector space.
We view $\Spin(\vb)$ as a subgroup of the group of invertible elements
of the Clifford algebra of $\vb$ (see e.g.\ \cite{lawson.michelsohn:89} or
\cite{hijazi:01}).
Let $e_1,\dots,e_n$
be an orthonormal basis of $\vb$. The complex spinor representation
$\Si$ of $\Spin(\vb)$ has dimension $2^{[n/2]}$.

\begin{proposition}
Let $g\in\Spin(\vb)$, $g\neq 1$. Then the multiplicity of the eigenvalue
$1$ of the endomorphism $g\in \End(\Si)$ is at most $2^{[n/2]-1}$.
\end{proposition}

\proof{}
We set $A_j:=e_{2j-1}\cdot e_{2j}\in Cl(\vb)$ for $j=1,\dots,m$, $m:=[n/2]$.

Any $g\in \Spin(\vb)$ is contained in a maximal torus,
i.e.\ there is an $h\in \Spin(\vb)$ and $t_j\in \mR$
such that
  $$g = h\cdot\exp(t_1A_1)\cdot\exp(t_2A_2)\cdot\dots \cdot\exp(t_mA_m)\cdot h^{-1}.$$
Let $h'$ be the image of $h$ under the map $\Spin(\vb)\to \SO(\vb)$.
Then $h\cdot\exp(t_jA_j)\cdot h^{-1}= \exp ( t_j h \cdot A_j \cdot h^{-1})$.
We set
  $$\widehat A_j:=h \cdot A_j \cdot h^{-1}= h'(e_{2j-1})\cdot h'(e_{2j})$$
and
  $$g_j:=\exp(t_j\widehat A_j)= \cos t_j + \sin t_j \widehat A_j.$$
The $\widehat A_j$ are pairwise commuting anti-self-adjoint endomorphisms.
Hence they are simultaneously
diagonalizable, with eigenvalues $i$ and $-i$.
Furthermore $h'(e_{2j})$ anti-commutes with $\widehat A_j$
and commutes with $\widehat A_k$, $j\neq k$.
Hence, all simultaneous eigenspaces have the same dimension, which is $1$.

We conclude that all $g_j$ are simultaneously diagonalizable with eigenvalues
$\exp(it_j)$ and $\exp(-it_j)$,
having $1$-dimensional simultaneous eigenspaces.
Thus $g$ has the eigenvalues
  $$e^{i(\pm t_1\pm t_2 \pm \dots \pm t_m)}$$
where the signs vary independently, each sign combination providing an
eigenspace of multiplicity $1$.
As a consequence, for any $g\neq 1$, the multiplicity of the eigenvalue
$1$ is at most $2^{m-1}$.
\qed

\begin{proposition}
The fixing dimension of the complex spinor representation of $\Spin(n)$
is
  $$r:=
     \begin{cases}
        1 & \mbox{if $n=2,3$}\cr
        2^{\left[{n\over 2}\right]-1}+1 & \mbox{if $n\geq 4$}
     \end{cases}$$
\end{proposition}
\proof{}If $n=2$, the spinor representation is
  $$S^1\to \SU(2), \quad z\mapsto
   \begin{pmatrix}z & 0 \cr 0 & z^{-1}\end{pmatrix}.$$
There are no invariant subspaces, hence $r=1$.

If $n=3$, then the spinor representation is the identity
$\Spin(3)=\SU(2)\to \SU(2)$. Let $W$ be a $1$-dimensional subspace of $\mC^2$.
Any $h\in \Stab^{\SU(2)}(W)$ can be diagonalized with
eigenvalues $\la$ and $\la^{-1}$. However, as $W$ is fixed by $h$, we obtain
$\lambda=1$, and hence $h=1$. We have thus shown that $r\leq 1$. Obviously
$r\geq 1$.

If $n\geq 4$, then example (\ref{kthreeex}) in the introduction
shows that $r>2^{\left[{n\over 2}\right]-1}$.
Together with the previous proposition we obtain the
result.
\qed

%%%%%%%%%%%%%%%%%%%%%%%%%%%%%
\section{Proof of Theorem~\ref{theo.oneone}}
%%%%%%%%%%%%%%%%%%%%%%%%%%%%%

We briefly recall some definitions from spin geometry.
Details can be found for example in
\cite{lawson.michelsohn:89}.

\begin{definition}
Let $(M,g)$ be a Riemannian manifold. Let $\PSO(M)$ be the frame bundle
over $M$. A \emph{spin structure} is a
$\Spin(n)$-principal bundle $\Pspin(M)$ together with a
$\Theta:\Spin(n)\to SO(n)$-equivariant fiber map
$\chi:\Pspin(M)\to \PSO(M)$ over the identity $M\to M$.
\end{definition}

\begin{example}
Let $G$ be an $n$-dimensional Lie group and $\Gamma$ a lattice in $G$.
The frame bundle of $G$ is trivialized by left invariant frames, i.e.\
$\SO(G)=G\times \SO(n)$. Hence, there is a spin structure on $G$
given by $\Spin(G)=G \times \Spin(n)$ where $\chi$ is the identity in the
first component and the standard map $\Spin(n)\to \SO(n)$ in the second
component. The frame bundle
of $\Gamma\backslash G$ is $(\Gamma\backslash G)\times\SO(n)$.
One possible spin structure
on $\Gamma\backslash G$ is $(\Gamma\backslash G)\times\Spin(n)$ together
with the equivariant map $\id \times \Theta$.
This spin structure is called the \emph{trivial spin structure}.
%Note that the term ``trivial'' used
%here is different from the fact that $(M,\chi)$ is spin-cobordant zero.
\end{example}

\begin{definition}
Let $\rho:\Spin(n)\to U(\Si)$ be the complex spinor representation.
The \emph{spinor bundle} is defined as the associated vector bundle
  $$\Si M:=\Pspin(M)\times_\rho \Si.$$
\end{definition}

As $\rho$ is a $2^{[n/2]}$-dimensional complex representation,
the complex vector bundle $\Si$ has rank $k:=2^{[n/2]}$.
The holonomy is contained in $\Spin(n)$,
the inclusion $\Spin(n)\hookrightarrow \U(k)$ given by the
spinor representation.

As before we denote the fixing dimension of the spinor
representation by
  $$r=r(n)=
     \begin{cases}
        1 & \mbox{if $n=2,3$}\cr
        2^{\left[{n\over 2}\right]-1}+1 & \mbox{if $n\geq 4$}.
     \end{cases}$$

\begin{theooneone}
\contentoneone
\end{theooneone}

\proof{}
We will use the Abresch's smoothing theorem
\cite[Theorem~1.12]{cheeger.fukaya.gromov:92}
For the convenience of the reader we will summarize
the smoothing theorem in Appendix~\ref{appendix.smoothing}.
This theorem states that for any $\de>0$ there is a constant
and $K_1(K,n,\de)$ such that any metric $g$ on a
compact manifold $M^n$ with $|\sec_g|<K$
can be approximated by another metric $\tilde g$ on $M$ with
\begin{enumerate}[{\rm (1)}]
\item $e^{-\de}g\leq \ti g \leq e^{\de}g$,
\item $|\na^g -\na^{\ti g}|\leq \de$,
\item $|\sec_{\ti g}|\leq K+\delta$, and
\item $|\na R_{\ti g}|\leq K_1$.
\end{enumerate}
As the eigenvalues of the connection Laplacian $\na^*\na$ on the spinor bundle
are uniformly continuous under $C^1$-perturbations of the metric
(Proposition~\ref{prop.conn.cont}), this shows that it is sufficient to prove the
theorem under the additional assumption that $|\nabla R|$ is
bounded. We will formulate the remaining step as a
lemma.
\qed

\begin{lemma}
Let $(M^n,g, \chi)$ be a compact Riemannian spin manifold with
$|\Sec|<K$, $|\nabla R|<K$, $\diam < D$. Then there
is $\varepsilon=\varepsilon(n,K,D)>0$,
such that if  $\lambda_r(\nabla^*\nabla)<\varepsilon$,
then $M$ is diffeomorphic to a nilmanifold.
Furthermore, $\chi$ is the trivial spin structure on $M$.
%There is an $\ep= \ep(K_0,K_1,D,n)>0$
%such that any compact Riemannian spin manifold
%$(M^n,g,\chi)$ which satisfies
%\begin{enumerate}[{\rm (i)}]
%\item $|\sec|< K_0$ and $\diam<D$,
%\item $|\na R|< K_1$
%\item the connection-Laplacian $\nabla^*\nabla$ on the spinor bundle
%has at least $r(n)$ eigenvalues smaller than $\ep$.
%\end{enumerate}
%is diffeomorphic to a nilmanifold and the $\chi$
%is the trivial spin structure.
\end{lemma}

\proof{of the lemma}
We apply Proposition~\ref{prop.hol.absch} for $\vb=\Si M$, $H=\Spin(n)$.
We obtain a frame $\cE$ of $\Si M$ with $|\na \cE|=\tau(\lambda|n,K,D)$.
The spin structure $\chi:\Pspin (M)\to \PSO(M)$ maps $\cE$ to $\chi(\cE)$
with
  $$|\na\chi(\cE)|\leq \tau(\epsilon|n,K,D),$$
i.e.\ an almost parallel frame of $TM$.
Now, using \cite{ghanaat:89}
we see that $M$ is $C^0$-close and diffeomorphic to a nilmanifold $\Gamma \backslash N$ with $\Gamma$ a
cocompact lattice in the nilpotent Lie group $N$.
Let ${\cF'}$ be a frame on $M$ which is sufficiently close to $\chi(\cE)$ and
whose pullback to $N$ is left-invariant.
It can be lifted, i.e.\ there is
a frame $\cE'$ with $\chi(\cE')=\cF'$, hence the spin structure is trivial.
\qed

%%%%%%%%%%%%%%%%%%%%%%%%%%%%%%%%%%%%%%%%%%%%%%
\appendix
%%%%%%%%%%%%%%%%%%%%%%%%%%%%%%%%%%%%%%%%%%%%%%
\def\secvar{S}
\section{Some analytical tools}\label{appendix.ana.tools}

Here we outline the results from
\cite{petersen.sprouse:99},\cite{petersen.sprouse:erratum} which we need.
The main analytic tool is the following lemma which follows from
Moser iteration. Note that Lemma 3.1 in \cite{petersen.sprouse:99}
is incorrect. A correct version is as follows (\cite{petersen.sprouse:erratum}).
Similar bounds were also obtained in \cite{aubry.colbois.ghanaat.ruh:03} and
\cite{ballmann.bruening.carron:02},
where a version of Lemma A.4 was derived which does not depend on $|\divv R^V|$.

%\begin{lemma}\label{lemma.one.one}
%Let $(M,g)$ satisfy $\Ric \geq -k^2$, $\diam < D$. Then
%for a funcion $u$ on $M$ satisfying $\Delta u \leq \alpha u$,
%$\alpha\geq 0$
%we have $\|u\|_\infty \leq (1 +\tau(\alpha|n,k,D))\|u\|_2$.
%\end{lemma}

\begin{lemma}\label{lemma.one.one}
Let $(M,g)$ satisfy $\Ric \geq -k^2$, $\diam < D$. Then
for a funcion $u$ on $M$ satisfying $\Delta u \leq \alpha u + \beta$,
$\alpha,\beta\geq 0$
we have $\|u\|_\infty \leq \tau(\|u\|_2|\alpha,\beta,n,k,D)$.
If $\beta=0$ then in fact
$\|u\|_\infty \leq (1 + \tau(\alpha|n,k,D))\|u\|_2$.
\end{lemma}

Here the diameter and Ricci curvature bounds give a bound on the Sobolev
constant used in Moser iteration by a result of Gallot, and the lower Ricci curvature
bound is implied the bounds on sectional curvature which we have assumed.
%\begin{lemma}
%Let $(M,g)$ satisfy $\Ric \geq -K^2$, $\diam < D$.
%Suppose $u$ is a function on $M$ with $\Delta u = f$. Then
%$\| u - \overline u \|_\infty \leq C(n,K,D)\| f \|_2$.
%\end{lemma}
Then, a standard argument yields the following.
\begin{lemma}\label{lematwo}
Let $V$ be a vector bundle over $M$. Suppose $M$ has
$\Ric \geq -k^2$, $\diam < D$. Then
for any section of $V$ satisfying
$\langle \nabla^*\nabla
\secvar, \secvar \rangle \leq \lambda |\secvar|^2$, $\|\secvar\|_2 =
1$,
we have $\|\secvar\|_\infty \leq 1 + \tau(\lambda|n,k,D)$.
\end{lemma}
\proof{}
\begin{eqnarray*}
\Delta |\secvar|^2 & = &2\langle \nabla^*\nabla \secvar , \secvar\rangle -
2|\nabla \secvar|^2\\
& \leq & 2\langle \nabla^*\nabla \secvar , \secvar
\rangle-2\Bigl|\nabla|\secvar|\Bigr|^2 \nonumber\\
& \leq & 2\lambda|\secvar|^2-2\Bigl|\nabla|\secvar|\Bigr|^2\nonumber
\end{eqnarray*}
Now we can use that we also have
\begin{equation*}
\Delta|\secvar|^2 = 2|\secvar|\;\Delta|\secvar| - 2\Bigl|\nabla|\secvar|\Bigr|^2
\end{equation*}
and solving for $\Delta |\secvar|$ we get $\Delta |\secvar| \leq \lambda |\secvar|$,
and hence \ref{lemma.one.one} gives the desired result.
\qed
To bound $|\nabla \secvar|$ we apply the following Bochner formula
\begin{lemma}\label{lemma.bochner}
Let $R^V:TM\otimes TM\otimes V\to V$ denote the curvature of the vector
bundle $V$. Let $\divv^1 R^V$ be minus
the metric contraction of $\nabla R^V$
in the first two slots.
Then for any section $\secvar$ of $V$
$$ \nabla(\nabla^*\nabla) \secvar =
    (\nabla^*\nabla)\nabla \secvar  - (\divv^1 R^V) \secvar
    + \nabla_{\Ric(.)}\secvar+ 2 c_{12}(\id \otimes R^V)(\na\secvar)$$
\end{lemma}

\proof{}
First note that, as the metric is parallel, metric contraction is parallel.
The metric contraction of the $i$-th slot with the $j$-th slot is denoted by
$c_{ij}$.
Let $\tau:TM\otimes TM\to TM\otimes TM, \quad X\otimes Y \mapsto Y\otimes X$.
Note that $R^V \secvar= \na\na \secvar - (\tau\otimes \id)\na\na \secvar$.
We calculate
\begin{eqnarray*}
  \na(\na^*\na)\secvar & = & -\nabla c_{12} (\na \na \secvar) \\
                    & = & -c_{23} (\na\na\na \secvar)\\
                    & = & -c_{23} (\tau\otimes \id \otimes \id)
                    (\na\na\na \secvar+ R^{T^*M \otimes V}\nabla\secvar)\\
                    & = & -c_{23} R^{T^*M \otimes V}\na \secvar -c_{13}\na\na\na\secvar
\end{eqnarray*}
The first summand gives
   $$-c_{23} R^{T^*M \otimes V}\na \secvar = c_{23} \na_{R(\cdot,\cdot)\cdot}\secvar
     -c_{13} (\id \otimes R^V)(\na \secvar) = \na_{\Ric(\cdot)} \secvar + c_{12} (\id \otimes R^V)(\na \secvar).$$
For the second term,
\begin{eqnarray*}
-c_{13}\na\na\na\secvar &=&-c_{13}\na\left((\tau\otimes \id)\na\na\secvar+ R^V\secvar\right)\\
& = & -c_{13}(\id\otimes \tau\otimes \id)(\na\na\na \secvar) + c_{12}\na R^V\secvar\\
& = & (\na^*\na)\na\secvar -\divv^1 R^V \secvar  +c_{12}(\id \otimes R^V)(\na \secvar),
\end{eqnarray*}
where we have used the definition $\divv^1 R^V\secvar= -c_{12}(\na R^V)\secvar$.
\qed
As a consequence of the Lemmas~\ref{lemma.one.one} and \ref{lemma.bochner}
we have
\begin{lemma}\label{lemafour}
Suppose that $M, V$ have $\Ric \geq -k^2$,
$|c_{12}(\id \otimes R^V)|,|\divv R^V| < K$, and
$diam < D$. Let $\secvar$ be an eigensection of $V$ with
$\nabla^*\nabla \secvar = \lambda \secvar$, and
$\|\secvar\|_2=1$.
Then $\| \nabla \secvar\|_\infty \leq \tau(\lambda | n,k,K,D)$.
\end{lemma}
Note that $|\Ric|$ and $|c_{12}(\id \otimes R^V)|$ are bounded by $|\sec|$
and $|R^V|$, and $|\divv R^V|$ is
bounded by $|\nabla R^V|$.
\proof{}
First of all
\begin{eqnarray*}
\Delta | \nabla \secvar |^2 &=&
2\langle \nabla^*\nabla(\nabla \secvar) , \nabla \secvar \rangle - 2|\nabla
\nabla \secvar|^2\\
&\leq&2\langle \nabla^*\nabla(\nabla \secvar) , \nabla \secvar \rangle -
2\Bigl|\nabla|\nabla \secvar|\Bigr|^2\nonumber
\end{eqnarray*}
Again we can use
\begin{equation}
\Delta |\nabla \secvar|^2 = 2|\nabla \secvar|\;\Delta|\nabla\secvar| -
2\Bigl|\nabla|\nabla\secvar|\Bigr|^2
\end{equation}
Which gives
\begin{eqnarray*}
|\nabla \secvar|\Delta |\nabla \secvar | &\leq &\langle \nabla^*\nabla (\nabla \secvar), \nabla \secvar \rangle\\
& = & \langle \nabla(\nabla^*\nabla) \secvar - \nabla_{\Ric(.)}\secvar + \divv R^V \secvar  - 2c_{12}(\id \otimes R^V)(\na\secvar), \nabla \secvar\rangle\nonumber\\
& \leq & (\lambda + k^2+2|c_{12}(\id \otimes R^V)|)|\nabla \secvar|^2 + |\divv R^V||\secvar||\nabla\secvar|.\nonumber
\end{eqnarray*}
Hence
\begin{eqnarray*}
\Delta |\nabla \secvar | & \leq & (\lambda + k^2+2|c_{12}(\id \otimes R^V)|)|\nabla \secvar| + |\divv R^V||\secvar|\\
&\leq& (\lambda + k^2 + 2|c_{12}(\id \otimes R^V)|)|\nabla \secvar| + |\divv R^V|(1 + \tau(\lambda|n,k,D)).\nonumber
\end{eqnarray*}
Then finally we can use Lemma~\ref{lemma.one.one}, along with the fact that
$\| \nabla \secvar \|_2 = \lambda\|\secvar\|_2$.
\qed
We then note that Lemma~\ref{lemma.bochner} gives us a bound on
$\nabla^*\nabla (\nabla \secvar)$ from which we can conclude from Lemma A.4
that $\int_M \langle \nabla^*\nabla (\nabla \secvar), \nabla \secvar \rangle dV =
\| \nabla \nabla \secvar_i\|_2$ is small.
%v2.3
Hence, Theorem~\ref{theo.S.est} is proven.

Finally we include a proof of Theorem~\ref{theo.S.orth}.

\begin{theosorth}
Suppose that $S_1,\dots,S_m$ are $L^2$-orthonormal eigensections of $\vb$, with
eigenvalues $\lambda_1 \leq \dots\leq \lambda_m$.
Then with for $k,K,D$ as above
$$ \| \langle S_i, S_j \rangle - \delta_{ij} \|_\infty \leq \tau(\lambda_m|n,K,k,D).$$
%\contenttheosorth
\end{theosorth}

\proof{}For any unit length vector $X$ we calculate for
\begin{eqnarray*}
 |\na_X\<S_i,S_j\>|&=&|\<\na_X S_i, S_j\> + \<S_i,\na_X S_j\>|\\
  &\leq & \|\na S_i\|_\infty \|S_j\|_\infty + \|S_i\|_\infty \|\na
  S_j\|_\infty\\
  & \leq & \tau(\max(\la_i,\la_j)|n,k,K,D)
\end{eqnarray*}
where we used Lemmata~\ref{lematwo} and \ref{lemafour} in the last
inequality.
Together with $\int_M \<S_i,S_j\>=\delta_{ij}$ the statement
easily follows.
\qed

\section{Connection Laplacians under perturbations of the metric}\label{appendix.pertcl}
We assume here that $M$ is a compact spin manifold with a fixed (topological)
spin structure. The topological spin structure defines for any metric
$g$ on $M$ a (metric) spin structure $\Pspin(M,g)\to \PSO(M,g)$.

\begin{proposition}\label{prop.conn.cont}
Let $\nabla^g$ be the Levi-Civita-connection on the spinor bundle with respect to the metric $g$.
Let $\la_1(g)\leq \la_2(g)\dots$ be
the eigenvalues of the connection Laplacian ${\nabla^g}^*{\nabla^g}$.
For two metrics $g$ and $\ti g$ let $\de=\de(g,\ti g)$ be the smallest number such that
  $$e^{-\de}\ti g(X,X)\leq g(X,X) \leq e^{\de}\ti g(X,X)\qquad
    \forall X\in TM,$$
  $$|\na^g-\na^{\ti g}|_g\leq \de.$$
Then
  $$e^{-{1001\over 1000}\,\de} \la_k(\ti g) - \tau(\de)\leq \la_k(g)\leq
     e^{{1001\over 1000}\, \de} \la_k(\ti g) + \tau(\de),$$
where $\tau(\delta)\to 0$ for $\delta \to 0$.
\end{proposition}

\proof{}
Let $A$ be the unique positive selfadjoint endomorphism of $TM$ such that
  $$\ti g(AX,AY)= g(X,Y).$$
As a consequence
  $$e^{-\de/2}\Id \leq A \leq e^{\de/2}\Id$$
in the sense of symmetric operators and
  $$|\witi\na(A^2)|_{\ti g}=|\witi\na  g|_{\ti g}
=|(\witi\na -\na) g|_{\ti g}\leq \tau_1(\de).$$
Hence also $|\witi \na A|_{\ti g} \leq \tau_2(\de)$.
Let $e_1,\ldots, e_n$ be a local orthonormal frame for $g$.
Then $Ae_1,\ldots, Ae_n$ is an orthonormal frame for $\ti g$.
The connection-1-forms $\om$ and $\ti\om$ are defined as
  $$\om(X)_j^k:= g(\na_X e_j,e_k)\qquad \witi\om(X)_j^k:= \ti g(\witi\na_X Ae_j,Ae_k).$$
We calculate
\begin{eqnarray*}
    |\om(X)_j^k -\witi\om(X)_j^k|
    & \leq & |g((\na_X-\witi\na_X) e_j,e_k)|+ |g(A^{-1}(\witi\na_X A)e_j,e_k)|\\
    & \leq & \left(|\na-\witi\na|_g + e^{\de/2}|\witi \na A|_{\ti g}\right)|X|
      \leq \tau_3(\de)\,|X|
\end{eqnarray*}
The map $A$ induces an $\SO(n)$ equivariant fiber map
$\PSO(M,g)\to \PSO(M,\witi g)$ which lifts to the spin structure
(if the spin structures coincide as topological spin structures).
The corresponding vector bundle map of the associated bundles is an
isomorphism of vector bundle $A:\Si(M,g)\to\Si(M,\witi g)$
which preserves length fiberwise and such that $A(V\cdot \phi)= A(V)\cdot A(\phi)$ for all $V\in T_pM$, $\phi\in \Si_p M$.
  $$A(\na_X\phi)-\witi\na_X (A \phi)= {1\over 4}\,
    \sum_{j,k=1}^n(\om(X)_j^k-\witi\om(X)_j^k)A(e_j\cdot e_k \cdot \phi).$$
Hence $|A(\na_X\phi)-\witi\na_X A \phi|\leq {n^2\tau_3(\de)\over 4}$.
  $$|d(\det A)|_{\witi g}\leq \tau(\de).$$
We set $\witi\psi:=(\det A) A\psi$. The map
$L^2(\Si(M,g))\to L^2(\Si(M,\ti g))$, $\psi\mapsto \witi \psi$ is an isometry.
Then
$$\<\na^*\na\psi,\psi\>^{1/2}\leq e^{\de/2}\<\witi\na^*\witi\na\witi\psi,\witi\psi\>^{1/2} + \tau(\de).$$

From this we deduce
  $$(1-\tau(\de))\,e^{-\de} \la_k(\ti g) - \tau(\de)\leq \la_k(g)\leq (1+\tau(\de))\,e^{\de} \la_k(\ti g) +  \tau(\de).$$
As a consequence
  $$e^{-{1001\over 1000}\,\de} \la_k(\ti g) - \tau(\de)\leq \la_k(g)\leq
     e^{{1001\over 1000}\, \de} \la_k(\ti g) + \tau(\de).$$
\qed

The proof runs completely analogous using that
  $$|A(X\cdot \phi)-X\cdot A(\phi)|
    = |((A-\Id)X)\cdot \phi|\leq \tau(\de)|X|\,|\phi|.$$

\section{Smoothing Riemannian metrics}\label{appendix.smoothing}
For the convenience of the reader we want to recall
the smoothing result that we need, which follows from work in
\cite{abresch:88}, \cite{bando:87} \cite{BMR:84}, \cite{cheeger.gromov:85}.
A brief survey on such results is contained in Section 5 of \cite{fukaya:90}.
In particular, Proposition 5.9 applied to Theorem 5.1 in \cite{fukaya:90}
gives:

%the smoothing result by Abresch
%\cite{abresch:88},
%Cheeger--Fukaya--Gromov \cite{cheeger.fukaya.gromov:92},
%Rong \cite{rong:96} et al.\ which we apply in the article.

\begin{theorem}\label{theo.smoothing}
%[{Abresch \cite[Theorem~1.12]{cheeger.fukaya.gromov:92}}]
For any $\delta=\delta(K,n)>0$ there is
%$K_0=K_0(K,n)$
$K_1=K_1(n,K,\delta)$ such that on any $n$-dimensional
complete Riemannian manifold $(M,g)$ with $|\sec|\leq K$ there is
a Riemannian metric $\ti g$ on $M$ such that
  $$e^{-\delta}\ti g(X,X)\leq g(X,X) \leq e^{\delta}\ti g(X,X)\qquad
    \forall X\in TM,$$
  $$|\na^g-\na^{\ti g}|_g\leq \delta.$$
  $$|\sec_{\ti g}|_{\ti g}\leq K + \delta, \qquad |\nabla^{\ti g}
    R_{\ti g}|_{\ti g}\leq K_1.$$
\end{theorem}

%\section{A metric on the cylinder}\label{appendix.cyl}

%\bibliographystyle{amsbernd}
%\bibliography{literatur}

\def\pages#1{#1}

%\bibliographystyle{amsbernd}
%\bibliography{literatur}

\begin{thebibliography}{100}
\bibitem{abresch:88}U. Abresch, {\it \"Uber das Gl\"atten Riemannscher Metriken},
Ha\-bi\-li\-ta\-tions\-schrift, Rheinische Friedrich-Wilhelms-Universit\"at, Bonn (1988).
\bibitem{ammann:02}B. Ammann, {\it Dirac eigenvalue estimates on 2-tori},
to appear in J. Geom. Phys.
\bibitem{ammann.baer:98}B. Ammann, C. B\"ar, {\it The Dirac Operator on Nilmanifolds and Collapsing Circle Bundles}, Ann. Global Anal. Geom. {\bf 16} no. 3 (1998), \pages{221--253}
\bibitem{ammann.sprouse:detailssmallev1:1}B. Ammann, C. Sprouse, {\it
Details to Examples in ``Manifolds with small Dirac eigenvalues are nilmanifolds''},
downloadable on \wwwlink{http://www.berndammann.de/publications/smallev1}.
\bibitem{anderson:90}M. Anderson, Convergence and rigidity of
manifolds under Ricci curvature bounds, Invent. Math {\bf 102} (1990),
\pages{429--445}.
\bibitem{aubry.colbois.ghanaat.ruh:03}E.~Aubry, B.~Colbois, P.~Ghanaat and
E.~Ruh,
{\it Curvature, Harnack's
inequality, and a spectral characterization of nilmanifolds},
Ann. Golb. Anal. Geom. {\bf 23} (2003), \pages{227--246}
\bibitem{baer:92}
C. B\"ar,
{\it Lower eigenvalue etimates for Dirac operators},
Math.\ Ann.\ {\bf 293} (1992), \pages{39--46}
\bibitem{baer:93}C. B\"ar, {\it Real Killing Spinors and Holonomy},
Comm. Math. Phys. {\bf 154} (1993), \pages{509--521}
%\bibitem{baer:96a}C. B\"ar, The Dirac operator on space
%forms of positive curvature. J. Math. Soc. Japan {\bf 48}
%(1996) 69-83.
\bibitem{baer:96b}
C.~B{\"a}r, \emph{Metrics with harmonic spinors}, Geom. Funct. Anal. \textbf{6}
  (1996), \pages{899--942}
\bibitem{baer.dahl:02}
C.~B{\"a}r, M.~Dahl, \emph{{Surgery and the Spectrum of the Dirac
  Operator}}, J. reine angew. Math. \textbf{552} (2002), \pages{53--76}.
\bibitem{baer.dahl:p03}C.~B\"ar, M.~Dahl,
\emph{The first Dirac eigenvalue on manifolds with positive scalar curvature},
Preprint 2003, to appear in Proc.\ AMS.
\bibitem{ballmann.bruening.carron:02}
W.~Ballmann, J.~Br\"uning, G. Carron, \emph{Eigenvalues and holonomy},
International Math. Research Notices \textbf{12} (2003), \pages{657--665}
\bibitem{bando:87} S. Bando, \emph{{Real analyticity of solutions of Hamilton's
equation}}, Math. Z. \textbf{195} (1987), \pages{93--97}.
\bibitem{BMR:84}J. Bemelmans, M. Min-Oo, E. Ruh, Smoothing
Riemannian metrics, Math Z. {\bf 188} (1984), \pages{69--74}.
\bibitem{cheeger.gromov:85} J. Cheeger, M. Gromov, \emph{On the
characteristic numbers of complete manifolds of bounded curvature and
finite volume}, Differential Geometry And Complex Analysis, H. E. Rauch Memorial
Volume (1985), \pages{115--154}.
\bibitem{cheeger.fukaya.gromov:92}
J.~Cheeger, K.~Fukaya, M.~Gromov, \emph{{Nilpotent Structures and invariant
  metrics on collapsed manifolds.}}, J. Amer. Math. Soc. \textbf{5} (1992),
\pages{327--372}.
\bibitem{CC:90}B. Colbois, G. Courtois,
\emph{A note on the first nonzero eigenvalue of the Laplacian acting
on $p$-forms}, Manuscripta Math. {\bf 68} (1990), \pages{143--160}.
%\bibitem{dahl:p03}M. Dahl,
%\emph{Prescribing eigenvalues of the Dirac operator},
%Preprint, ArXiv math.DG/0311172.
\bibitem{friedrich:80} T. Friedrich, Der erste Eigenwert des Dirac-Operators einer
kompakten Riemannschen Mannigfaltigkeit nichtnegativer Skalarkr\"ummung,
Math. Nachr. {\bf 97} (1980), \pages{117--146}.
\bibitem{friedrich} T. Friedrich, Dirac Operators in Riemannian Geometry,
Graduate Studies in Mathematics {\bf 25} (2000).
\bibitem{fukaya:90} K. Fukaya, \emph{Hausdorff convergence of Riemannian manifolds
and its applications}, Recent Topics in Differential and Analytic Geometry, Adv.
Stud. Pure Math. \textbf{18-I} (1990), \pages{143--238}.
\bibitem{gao:90} L. Gao, \emph{Convergence of Riemannian manifolds:
Ricci and $L^\frac{n}{2}$ curvature pinching}, J. Diff. Geom. {\bf 32}
(1990), \pages{349--381}.
\bibitem{ghanaat:89} P. Ghanaat, \emph{Almost Lie groups of type $\mathbb R^n$},
J. Reine Angew. Math. \textbf{401} (1989), \pages{60--81}.
\bibitem{hijazi:01} O. Hijazi,
{\it Spectral Properties of the Dirac Operator and Geometrical Structures},
Proceedings of the Summer School on Geometric
Methods in Quantum Field Theory, Villa de Leyva, Colombia,  July 12--30,
(1999), World Scientific 2001.
\bibitem{hitchin:74}
N.~Hitchin, \emph{Harmonic spinors}, Adv. Math. \textbf{14} (1974), \pages{1--55}.
\bibitem{lawson.michelsohn:89}
H.-B. Lawson, M.-L. Michelsohn, \emph{Spin geometry}, Princeton University
  Press, Princeton, 1989.
\bibitem{lott:02}
J.~Lott, \emph{Collapsing and {D}irac type operators}, Geom. Dedic. {\bf 91}, (2002) \pages{175-196}.
\bibitem{maier:97}
S.~Maier, \emph{Generic metrics and connections on spin- and
  spin-$\,^c$-manifolds}, Comm. Math. Phys. \textbf{188} (1997), \pages{407--437}.
\bibitem{min-oo.ruh:90}M. Min-Oo, E.Ruh, \emph{$L^2$-curvature
pinching}, Comment. Math. Helv. {\bf 65} (1990), \pages{36--51}.
\bibitem{petersen:97}P. Petersen, \emph{Convergence
Theorems in Riemannian Geometry}, in \emph{Comparison Geometry} Grove,
Petersen (Ed.), MSRI Publications {\bf 30} (1997), 167--202.
\bibitem{petersen.sprouse:99}P.~Petersen, C.~Sprouse,
\emph{Eigenvalue pinching for Riemannian vector bundles},
J.~Reine Angew.\ Math. {\bf 511} (1999), 73--86.
\bibitem{petersen.sprouse:erratum}P.~Petersen, C.~Sprouse,
\emph{Erratum to ``Eigenvalue pinching for Riemannian vector bundles''},
Preprint 2003.
\bibitem{petersen.sprouse:01}P.~Petersen, C.~Sprouse,
\emph{Eigenvalue pinching on $p$-forms} in
Proceedings of the Fifth Pacific Rim Geometry Conference.
Tohoku Mathematical Publications {\bf 20} (2001), 139--145.
\bibitem{pfaeffle:00}F.~Pf\"affle,
The Dirac spectrum of Bieberbach manifolds, J. Geom. Phys. {\bf 35}
(2000), 367--385
%\bibitem{rong:96}
%X.~Rong, \emph{{On the fundamental groups of manifolds of positive sectional
%curvature.}}, Ann. Math., II. Ser. \textbf{143} (1996), 397--411
\end{thebibliography}
\end{document}